\documentclass{amsart}
\usepackage{ amsmath, amsthm, amsfonts, hyperref, graphicx, ifpdf}
\usepackage[dvipsnames,usenames]{color}
\hypersetup{
   unicode=false,          
   pdftoolbar=true,        
   pdfmenubar=true,        
   pdffitwindow=false,     
   pdfstartview={},    
   pdftitle={},    
   pdfauthor={},     
   pdfsubject={},   
   pdfcreator={},   
   pdfproducer={}, 
   pdfkeywords={}, 
   pdfnewwindow=true,      
   colorlinks=true,       
   linkcolor=blue,          
   citecolor=blue,        
   filecolor=blue,      
   urlcolor=blue          
}

\newcommand{\R}{{\mathbb R}}

\newtheorem{theorem}{Theorem}[section]
\newtheorem{lemma}[theorem]{Lemma}
\newtheorem{pro}[theorem]{Proposition}
\newtheorem{cor}[theorem]{Corollary}

\theoremstyle{definition}

\theoremstyle{remark}
\newtheorem{remark}[theorem]{Remark}

\numberwithin{equation}{section}

\newcommand{\ee}{evolution equation}
\newcommand{\mc}{mean curvature}
\newcommand{\mcf}{mean curvature flow}
\newcommand{\maxp}{maximum principle}

\newcommand{\sap}{surface area preserving}
\newcommand{\sapmcf}{surface area preserving mean curvature flow}
\newcommand{\sff}{second fundamental form}
\newcommand{\ub}{uniform bound}
\newcommand{\vpmcf}{volume preserving mean curvature flow}
\newcommand{\wrt}{with respect to}

\newcommand{\be}{\begin{equation}}
\newcommand{\ene}{\end{equation}}
\newcommand{\br}{\begin{remark}}
\newcommand{\er}{\end{remark}}
\newcommand{\bl}{\begin{lemma}}
\newcommand{\el}{\end{lemma}}
\newcommand{\bcor}{\begin{cor}}
\newcommand{\ecor}{\end{cor}}
\newcommand{\ben}{\begin{enumerate}}
\newcommand{\een}{\end{enumerate}}
\newcommand{\bp}{\begin{proof}}
\newcommand{\ep}{\end{proof}}
\newcommand{\bpo}{\begin{pro}}
\newcommand{\epo}{\end{pro}}
\newcommand{\beq}{\begin{equation*}}
\newcommand{\eeq}{\end{equation*}}
\newcommand{\bear}{\begin{eqnarray*}}
\newcommand{\eear}{\end{eqnarray*}}
\newcommand{\bt}{\begin{theorem}}
\newcommand{\et}{\end{theorem}}

\DeclareMathAlphabet{\mathcal}{OMS}{cmsy}{m}{n}

\numberwithin{equation}{section}

\allowdisplaybreaks

\def\XXint#1#2#3{{\setbox0=\hbox{$#1{#2#3}{\int}$}
    \vcenter{\hbox{$#2#3$}}\kern-.5\wd0}}

\makeatletter
\def\@citestyle{\m@th\upshape\mdseries}
\def\citeform#1{{\bfseries#1}}
\def\@cite#1#2{{%
  \@citestyle[\citeform{#1}\if@tempswa, #2\fi]}}
\@ifundefined{cite }{%
  \expandafter\leqt\csname cite \endcsname\cite
  \edef\cite{\@nx\protect\@xp\@nx\csname cite \endcsname}%
}{}
\makeatother

\begin{document}
\newcommand{\osc}{{\text{osc}}}
\newcommand{\Vol}{{\text{Vol}}}
\newcommand{\V}{{\text{V}}}
\newcommand{\nn}{{\text{n}}}
\newcommand{\cM}{{\cal{M}}}
\newcommand{\Ric}{{\text{Ric}}}
\newcommand{\RE}{{\text{Re }}}
\newcommand{\LL}{{\cal{L}}}
\newcommand{\diam}{{\text {diam}}}
\newcommand{\dist}{{\text {dist}}}
\newcommand{\Area}{{\text {Area}}}
\newcommand{\Length}{{\text {Length}}}
\newcommand{\Energy}{{\text {Energy}}}
\newcommand{\SSS}{{\bold S}}
\newcommand{\K}{{\text{K}}}
\newcommand{\Hess}{{\text {Hess}}}
\def\RR{{\bold  R}}
\def\SS{{\bold  S}}
\def\TT{{\bold  T}}
\def\CC{{\bold C }}
\newcommand{\dv}{{\text {div}}}
\newcommand{\kg}{{\text {k}}}
\newcommand{\expp}{{\text {exp}}}
\newcommand{\e}{{\text {e}}}
\newcommand{\eqr}[1]{(\ref{#1})}
\newcommand{\ec}{\varepsilon_c}
\newcommand{\rc}{\rho_c}
\newcommand{\sztp}{\Sigma^{0,2\pi}}

\parskip1ex






\title[Stability of the surface area preserving mean curvature flow]{Stability of the surface area preserving mean curvature flow in Euclidean space}

\author{Zheng Huang}
\address[Z. ~H.]{Department of Mathematics, The City University of New York, Staten Island, NY 10314, USA}
\address{The Graduate Center, The City University of New York, 365 Fifth Ave., New York, NY 10016, USA}
\email{zheng.huang@csi.cuny.edu}

\author{Longzhi Lin}
\address[L.~L.]{Department of Mathematics\\Rutgers University\\110 Frelinghuysen Road\\Piscataway, NJ 08854-8019\\USA}
\email{lzlin@math.rutgers.edu}

\date{\today}
\subjclass[2010]{Primary 53C44, Secondary 58J35}

\begin{abstract}
We show that the {\sapmcf} in Euclidean space exists for all time and converges exponentially to a round sphere, if initially the
$L^2$-norm of the traceless {\sff} is small (but the initial hypersurface is not necessarily convex).
\end{abstract}

\maketitle

\section{Introduction}
Let $M^n$ be a smooth, embedded, closed (compact, no boundary) $n$-dimensional manifold in $\R^{n+1}$, and we evolve
it by the {\sapmcf}, that is,
\be
\label{eq-mcf}
\frac{\partial F}{\partial t} = \left(1 - h H\right)\nu, \qquad F(\cdot,0) = F_0(\cdot)\,.
\end{equation}
Here $F_0: M^n \to \R^{n+1}$ is the initial embedding, and $H = H(x,t)$ is the {\mc} and $\nu = \nu (x,t)$ is the outward
unit normal vector of $M_t = F(\cdot, t)$ at point $(x,t)$ (for simplicity, we simply write $(x,t)\in M_t$). And the function $h$
is given by
\be\label{def-h}
h = h (t) = \frac{\int_{M_t} H \, d\mu}{\int_{M_t} H^2 \, d\mu} \,,
\ene
where $d\mu = d\mu_t$ denotes the surface area element of the evolving surface $M_t$ {\wrt} the induced metric $g(t)$.
Clearly we have $H\not\equiv 0$ on $M_0$ since there is no closed minimal hypersurface in Euclidean space. A good monotonicity property of the {\sapmcf} \eqref{eq-mcf} is that the surface area of $M_t$ remains unchanged and the volume of the $(n+1)$-dimensional region enclosed by $M_t$ is non-decreasing along the flow, see Corollary ~\ref{monotonicity}.

We denote $A = \{a_{ij}\}$ as the {\sff} of $M_t$ and its traceless part as {\AA} $= A - \frac{H}{n}g$. Then we have
$|\text{\AA}|^2 = |A|^2 - \frac{1}{n}H^2$. This quantity measures the roundness of the hypersurface.

In this paper, we prove the following theorem on the stability of this {\sapmcf}:
\bt\label{MainThm}
Let $M_t^n \subset \mathbb{R}^{n+1}, n \geq 2,$ be a smooth compact solution to the {\sapmcf} \eqref{eq-mcf} for
$t\in [0,T)$ with $T\leq \infty$. Assume that $h(0)>0$. There exists $\epsilon>0$, depending only on $M_0$ and $h(0)$,
such that if
\be\label{initialCondition}
\int_{M_0}|\text{\AA}|^2 \,d\mu \leq \epsilon \,,
\ene
then $T = \infty$ and the flow converges exponentially to a round sphere.
\et
\br
The idea of using an initial {\sff} condition to pursue convergence of the flow was probably first studied in the case of
Ricci flow (\cite{Ye93}), later in the case of K\"{a}hler-Ricci flow (\cite{Che06, CLW09} and others). The setting for the
stability of the {\vpmcf} was studied by Escher-Simonett (\cite{ES98}) and Li (\cite{Li09}), under different set of
conditions. Escher-Simonett's approach is a center manifold analysis while Li's approach is to apply a parabolic version
of the Moser iteration method. Our approach is similar to the idea of iteration in \cite{Ye93, Li09}, in the cases of Ricci
flow and {\vpmcf}, respectively. Howover, the analytical nature of our case, namely the {\sapmcf}, is much more
complicated than that of the {\vpmcf}, since the function $h(t)$ contains two integral terms both involving the {\mc}.
Our approach is expected to be able to use to investigate the more general mixed volume {\mcf} studied in (\cite{McC04}).
\er
\br
In \cite{McC03}, McCoy proved that the {\sapmcf} exists for all time and converges to a sphere if the initial hypersurface
is {\it strictly convex}. As in the case of {\vpmcf} initiated by Huisken in \cite{Hui87}, strict convexity of the initial
surface is essential. In our setting, we do not assume any convexity for the initial hypersurface. Under the conditions
of the Theorem ~\ref{MainThm}, evolving surfaces become {\it mean convex} instantly after flow starts (see equation
\eqref{meanconvex}).
\er
\subsection*{Outline of the proof:} Our strategy is conventional: based on the initial bounds, we prove bounds on some time
interval for several geometric quantities (Theorem ~\ref{key}), then we prove exponential decay for these quantities on
the time interval of the interest (Theorem ~\ref{exp1}), which allows us to obtain uniform bounds for these quantities on
the interval (Theorem ~\ref{exp2}), therefore we can repeat above arguments to extend the time interval (Theorem ~\ref{ext}),
and the amount of extension only depends on the initial conditions. Main theorem then follows.

\subsection*{Plan of the paper:} There are four sections. In \S 2, we collect {\ee}s for various geometric quantities
associated to this flow, and provide some classic results that will be used in the proof. The proof of the main theorem is
contained in the last two sections: we provide key estimates and prove exponential decay for $|\text{\AA}|$ and
other quantities in \S 3, and we use these estimates to prove the long-time existence and exponential convergence in \S 4.

\subsection*{Acknowledgements} The research of Z. H. is partially supported by a PSC-CUNY award, and
Provost's Research Scholarship of CUNY-CSI.
\section{preliminaries}
We collect some necessary preliminary results in this section. In \S 2.1, we obtain {\ee}s for some key quantities and
operators, many of which were derived in \cite{McC03}; in \S 2.2, we state and use Hamilton's interpolation inequalities
for tensors to obtain a $L^2$ estimate (Lemma ~\ref{interpolation01}) on the gradients of the tensor $\text{\AA}$. A version
of the parabolic {\maxp} is also stated.
\subsection{Evolution of geometric quantities}
We start with the short time existence of the {\sapmcf} \eqref{eq-mcf} that is guaranteed by a work of Pihan:
\bt (\cite{Pih98}) Let $M_0$ be a smooth embedded compact $n$-dimensional manifold in
$\mathbb{R}^{n+1}$. Assume that $H \ne 0$ at some point of $M_0$ and $h(0)>0$, then there exists
$T_0>0$ such that the {\sap} {\mcf} \eqref{eq-mcf} exists and is smooth for $t\in [0,T_0)$.
\et
We now collect and derive some evolution equations of several geometric quantities which will be used later.
These quantities are:
\ben
\item
the induced metric of the evolving surface $M_t$: $g(t) = \{g_{ij}(t)\}$;
\item
the {\sff} of $M_t$: $A(\bullet,t) = \{a_{ij}(\bullet,t)\}$, and its square norm given by
\beq
|A(\bullet,t)|^2 = g^{ij}g^{kl}a_{ik}a_{jl};
\eeq
\item
the {\mc} of $M_t$ {\wrt} the outward normal vector: $H(\bullet,t) = g^{ij}a_{ij}$;
\item
the traceless part of the {\sff}: {\AA} $= A - \frac{H}{n}g$;
\item
the surface area element of $M_t$: $d\mu_t = \sqrt{det(g_{ij})}$.
\een
\bl (\cite{McC03})
The metric of $M_t$ satisfies the evolution equation
\be\label{ev-g}
\frac{\partial}{\partial t} g_{ij} = 2(1-hH)a_{ij}\,.
\ene
Therefore,
\be
\frac{\partial}{\partial t} g^{ij} = -2(1-hH)a^{ij}
\ene
and
\be\label{mu}
\frac{\partial}{\partial t} (d\mu_t) = H(1-hH)d\mu_t.
\ene
Moreover, the outward unit normal $\nu$ to $M_t$ satisfies
\be
\frac{\partial \nu}{\partial t}  = h \nabla H\,.
\ene
\el
As an easy consequence of \eqref{mu}, we have
\bcor (\cite{McC03}) \label{monotonicity}
\ben
\item
The surface area $|M_t|$ of $M_t$ remains unchanged along the flow, i.e.,
\beq
\frac{d}{dt}\int_{M_t}\,d\mu = \int_{M_t}(1-hH)H\,d\mu = 0\,.
\eeq
\item
The volume of $E_t$, the $(n+1)$-dimensional region enclosed by $M_t$, is non-decreasing along the flow, i.e.,
\beq
\frac{d}{dt}\text{Vol}\,(E_t) = \int_{M_t} \,d\mu - \frac{\left(\int_{M_t} H\,d\mu\right)^2}{\int_{M_t} H^2\,d\mu} \ge 0\,.
\eeq
\een
\ecor

\br
In Euclidean space, among all closed hypersurfaces, the sphere is of the least surface area with fixed enclosed volume,
and as well as of the largest enclosed volume with fixed surface area. Therefore from this point of view, it is natural
to study the sphere via both the {\vpmcf} and the {\sapmcf}.
\er

\bt (\cite{McC03})
The {\sff} satisfies the following {\ee}:
\be
\frac{\partial}{\partial t} a_{ij} = h\Delta a_{ij} + (1-2hH)a_i^{\,m}a_{mj} + h |A|^2 a_{ij}\,,
\ene
where $a_i^{\,m} = g^{ml}a_{li}$.
\et

\bcor (\cite{McC03})\label{Evol-A-H}
We have the {\ee}s for $H$, $|A|^2$ and $|\text{\AA}|^2$:
\begin{itemize}
\item[(i)] $\frac{\partial}{\partial t} H = h\Delta H- (1-hH)|A|^2$;
\item[(ii)] $\frac{\partial}{\partial t} |A|^2 = h\left(\Delta |A|^2 - 2|\nabla A|^2 + 2|A|^4\right) - 2\text{tr}\left(A^3\right)$,
\end{itemize}
where $\text{tr}\left(A^3\right) = g^{ij}g^{kl}g^{mn}a_{ik}a_{lm}a_{nj}$\,. Therefore we also have
\begin{itemize}
\item[(iii)] $\frac{\partial}{\partial t} |\text{\AA}|^2 =
h\Delta |\text{\AA}|^2 - 2h|\nabla \text{\AA}|^2 +
2h|A|^2|\text{\AA}|^2-2\left(\text{tr}(\text{\AA}^3)+{\frac{2}{n}} H|\text{\AA}|^2\right)$,
where $|\nabla \text{\AA}|^2 = |\nabla A|^2 - \frac{1}{n}|\nabla H|^2$.
\end{itemize}
\ecor
\bp
The last equation here is equivalent to the one from \cite{McC03}. To see this, we used the following fact (see
page 335 of \cite{Li09}):
\beq
\text{tr}\left(A^3\right) -\frac1n |A|^2H = \text{tr}\left({\text{\AA}}^3\right) + \frac2n |\text{\AA}|^2H.
\eeq
\ep
We can then derive the {\ee}s for the square norm of the gradients of the {\sff}.
\bcor\label{hdA-1}
We have the evolution euqation for $|\nabla^m A|^2$:
\begin{align}\label{hdA-2}
\frac{\partial}{\partial t} |\nabla^m A|^2 = & h\Delta |\nabla^m A|^2 - 2h|\nabla^{m+1} A|^2 +
\sum_{i+j+k=m} \nabla^i A\ast \nabla^j A \ast \nabla^k A \ast \nabla^m A \notag\\
& + \sum_{r+s=m} \nabla^r A\ast \nabla^s A \ast \nabla^m A\,,
\end{align}
where $S\ast \Omega$ denotes any linear combination (involving $h$) of tensors formed by contraction on $S$
and $\Omega$ by the metric $g$.
\ecor
\bp
The time derivative of the Christoffel symbols $\Gamma_{jk}^i$ is equal to
\begin{align*}
\frac{\partial}{\partial t} \Gamma_{jk}^i &= \frac{1}{2}g^{il} \left\{\nabla_j\left(\frac{\partial}{\partial t} g_{kl}\right) +
\nabla_k\left(\frac{\partial}{\partial t} g_{jl}\right)-\nabla_l\left(\frac{\partial}{\partial t} g_{jk}\right)\right\}\\
&= g^{il} \left\{\nabla_j\left( (1- hH)a_{kl}\right) + \nabla_k\left((1- hH)a_{jl}\right)-\nabla_l\left((1- hH)a_{jk}\right)\right\}\\
&= A \ast_h \nabla A + \nabla A\,,
\end{align*}
where $\ast_h = \ast$ denotes the contraction on tensors involving $h$ in the coefficients. Here we have also used the {\ee}
for the metric, i.e., \eqref{ev-g}. Then we can proceed as in \cite[\S 13]{Ham82} (see also \cite[\S 7]{Hui84}) to get
\eqref{hdA-2}.
\ep

In addition, we will need the following lemma on the time-derivative of the function $h(t) =\frac{\int_{M_t} H \, d\mu}{\int_{M_t} H^2 d\mu}$:
\bl\label{h_t}
\beq
\frac{\partial}{\partial t} h = \frac{\int_{M_t}[-(1-2hH)(1-hH)|A|^2 + H^2(1-hH)^2 + 2h^2|\nabla H|^2\, ]d\mu}{\int_{M_t} H^2 \, d\mu}.
\eeq
\el
\bp
For the sake of completeness, we compute as follows:
\begin{align*}
&\frac{\partial}{\partial t} h =\frac{\partial}{\partial t} \left( \frac{\int_{M_t} H \, d\mu}{\int_{M_t} H^2 \, d\mu}\right)\\
&= \left(\int_{M_t} H^2 \, d\mu\right)^{-1}\left[ \int_{M_t} - (1-hH)|A|^2+ H^2(1-hH)\, d\mu\right]\\
&\ \ \ -\left(\int_{M_t} H^2 \, d\mu\right)^{-1}\left[\int_{M_t} -2h^2|\nabla H|^2- 2hH(1-hH)|A|^2 + hH^3(1-hH)\, d\mu\right]\\
&=  \frac{\int_{M_t}[-(1-2hH)(1-hH)|A|^2 + H^2(1-hH)^2 + 2h^2|\nabla H|^2\, ]d\mu}{\int_{M_t} H^2 \, d\mu}\,.
\end{align*}
\ep

\subsection{Interpolation inequalities and {\maxp}}
We will need the following Hamilton's interpolation inequalities for tensors.
\bt(\cite{Ham82})\label{Hamiltonlemma2}
Let $M$ be an $n$-dimensional compact Riemannian manifold and $\Omega$ be any tensor on $M$. Suppose
\beq
\frac{1}{p} + \frac{1}{q} = \frac{1}{r} \quad \text{with } r\ge 1\,.
\eeq
We have the estimate
\beq
\left(\int_{M} |\nabla \Omega|^{2r}\, d\mu\right)^{1/r} \leq (2r-2+n)\,
\left(\int_{M} |\nabla^2 \Omega|^p\, d\mu\right)^{1/p}\left(\int_{M} |\Omega|^{q}\, d\mu\right)^{1/q}\,.
\eeq
\et

\bt(\cite{Ham82})\label{Hamiltonlemma}
Let $M$ and $\Omega$ be the same as the Theorem \ref{Hamiltonlemma2}. If $1\leq i\leq n-1$ and $m \ge 0$,
then there exists a constant $C = C(n,m)$ which is independent of the metric and connection on $M$, such that
the following estimate holds:
\beq
\int_{M} |\nabla^i \Omega|^{2m/i}\, d\mu \leq C\, \max_{M} |\Omega|^{2(m/i-1)}\int_{M}|\nabla^m \Omega|^2\, d\mu\,.
\eeq
\et
As an application of these inequalities, we provide an estimate that will be used later.
\bl\label{interpolation01}
For any $m\geq 1$ we have the estimate
\begin{align*}
\frac{d}{dt}\int_{M_t} |\nabla^m A|^2 \,d\mu &+ 2h \int_{M_t} |\nabla^{m+1} A|^2\,d\mu
\leq C\,\max_{M_t}\left(|A|^2 + |A|\right)\int_{M_t} |\nabla^m A|^2 \,d\mu\,,
\end{align*}
where $C=C(n,m, |h|)$.
\el
\bp
By integrating \eqref{hdA-2} of Corollary \ref{hdA-1} and using the generalized H\"{o}lder inequality we have
\begin{align*}
&\frac{d}{dt}\int_{M_t} |\nabla^m A|^2 \,d\mu -\int_{M_t}(1-hH)H|\nabla^m A|^2 \,d\mu + 2h \int_{M_t} |\nabla^{m+1} A|^2\,d\mu \\
\leq &C\,\Bigg\{\left(\int_{M_t} |\nabla^i A|^{2m/i} \,d\mu\right)^{i/2m}\left(\int_{M_t} |\nabla^j A|^{2m/j} \,d\mu\right)^{j/2m} \left(\int_{M_t} |\nabla^k A|^{2m/k} \,d\mu\right)^{k/2m}\\
&+ \left(\int_{M_t} |\nabla^r A|^{2m/r} \,d\mu\right)^{r/2m}\left(\int_{M_t} |\nabla^s A|^{2m/s} \,d\mu\right)^{s/2m}\Bigg\}\left(\int_{M_t} |\nabla^m A|^2 \,d\mu\right)^{1/2}\,,
\end{align*}
with $i+j+k = r+s = m$.

Applying Lemma \ref{Hamiltonlemma} for tensor $A$, we get
\beq
\left(\int_{M_t} |\nabla^q A|^{2m/q} \,d\mu\right)^{q/2m} \leq C\, \max_{M_t}|A|^{1-q/m}\left(\int_{M_t} |\nabla^m A|^2 \,d\mu\right)^{1/2m}\,,
\eeq
where $q= i, j, k,r,s$.

Also note that
\bear
\int_{M_t}|(1-hH)H|\nabla^m A|^2 d\mu &\leq& \max_{M_t}\{|H| + |h|H^2\}\int_{M_t} |\nabla^m A|^2 d\mu   \\
& \leq& C(n, |h|)\max_{M_t}\left(|A|^2 + |A|\right)\int_{M_t} |\nabla^m A|^2d\mu.
\eear
Combining these inequalities thus completes the proof.
\ep

We will need the following version of the {\maxp}, especially in the proof of the Theorem ~\ref{key}.
\bt\label{maxPrin} (Maximum principle, see e.g. \cite[Lemma 2.12]{CLN06}) Suppose $u: M \times [0,T] \to \R$ satisfies
\beq
\frac{\partial}{\partial t} u  \leq a^{ij}(t)\nabla_i\nabla_j u + \langle B(t), \nabla u\rangle + F(u)\,,
\eeq
where the coefficient matrix $\left(a^{ij}(t)\right)>0$ for all $t\in [0,T]$, $B(t)$ is a time-dependent vector field and $F$ is a
Lipschitz function. If $u\leq c$ at $t=0$ for some $c>0$, then $u(x,t)\leq U(t)$ for all $(x,t)\in M_t, t\ge 0$, where $U(t)$
is the solution to the following initial value problem:
\beq
\frac{d}{dt} U(t) = F(U) \quad \text{with} \quad U(0) = c\,.
\eeq
\et

\section{Proof of Theorem \ref{MainThm}: estimates}

We break our proof into two sections. In this section, we provide key estimates that will be needed: in \S 3.1, we establish
the $L^\infty$-bound for {\AA} from its $L^2$ bound; in \S 3.2, we prove the exponential decay for $|\text{\AA}|$.
\subsection{Establishing bounds for geometric quantities}

Let us start with a result of Topping which plays an important role in the key estimates we will focus on in this subsection.
\bl(\cite{Top08})\label{top}
Let $M$ be an $n$-dimensional closed, connected manifold smoothly immersed in
$\mathbb{R}^N$, where $N \geq n+1$. Then the intrinsic diameter and the {\mc} $H$ of $M$ are related by
\beq
\text{diam}\,(M) \leq C(n)\int_M |H|^{n-1}\,d\mu\,.
\eeq
\el
We now prove the following key estimates. This allows us to obtain the $L^\infty$-bound for {\AA}, $|\nabla H|$
and $|1- hH|$ on some time interval. More specifically,

\bt\label{key}
Let $M_t^n \subset \mathbb{R}^{n+1}, n\geq 2,$ be a smooth compact solution to the {\sapmcf} \eqref{eq-mcf} for $t\in [0,T)$
with $T\leq \infty$. Assume that
\be\label{Condition-1}
 \max_{M_0}|A| \leq  \Lambda_0 \quad \text{and}\quad \frac{1}{\Lambda_0} \leq
 \left\{ h(0)\,, \,\int_{M_0}H^2\,d\mu\,,\, \int_{M_0}|\nabla^m A|^2\,d\mu \right\}\leq  \Lambda_0,
\ene
for some $\Lambda_0 \geq 2$ sufficiently large and all $m \in [1, \widehat{m}]$ with $\widehat{m}$ sufficiently large. Then there exists
$\epsilon_0 = \epsilon_0(n, |M_0|, \Lambda_0)>0$ and $T_1 = T_1(\Lambda_0) \leq 1$, such that if
\be\label{Condition-2}
\int_{M_0}|\text{\AA}|^2 \,d\mu \leq \epsilon \leq \epsilon_0\,,
\ene
then for all $t\in [0,T_1]$ we have
\be\label{estimate-1}
\max_{M_t}|A| \leq  2\Lambda_0\quad \text{and}\quad \frac{1}{2\Lambda_0} \leq \left\{ h(t)\,, \,\int_{M_t}H^2\,d\mu  \right\}\leq  2\Lambda_0\,.
\ene
Moreover, there exists $C_1 = C_1 (n, |M_0|, \Lambda_0)$ and some universal constant $\alpha \in (0,1)$ such that for any $t\in [0,T_1]$
\be\label{estimate-2}
\max_{M_t}\left(|\text{\AA}|+ |\nabla H| + |1- hH|\right) \leq C_1 \epsilon^{\alpha}\,.
\ene
\et

\br
It is very important to keep track of the dependence of constants on geometric quantities. As we shall see from the proof below,
the constant $C_1 = C_1 (n, |M_0|, \Lambda_0)$ is non-decreasing in $\Lambda_0$.
\er

\bp
By the short time continuity, we first let $t_1>0$  be the maximal time such that for all $t\in [0,t_1]$ we have
\be\label{Hsqbound}
\max_{M_t}|A| \leq 2\Lambda_0 \quad\text{and} \quad\frac{1}{2\Lambda_0} \leq \left\{ h(t)\,, \,\int_{M_t}H^2\,d\mu  \right\}\leq  2\Lambda_0\,.
\ene
Now using the fact that $|\text{tr}\left(A^3\right)|\leq |A|^3$ (see Lemma 2.2 \cite{HS99}), and Kato's inequality
$|\nabla |A||\leq |\nabla A|$, we derive from (ii) of Corollary \ref{Evol-A-H} to find
\beq
\frac{\partial}{\partial t} |A| \leq  h \Delta |A| + 2\Lambda_0 |A|^3 + |A|^2 \quad \text{on }\, M_t \text{ for all } t\in [0,t_1]\,.
\eeq
Then by the {\maxp} (Theorem \ref{maxPrin}), we have:
\beq
\max_{M_t}|A| \leq U(t) \text{ for all } t\in [0,t_1]\,,\text{ with } U(0) = \Lambda_0,
\eeq
where $U(t) >0$ solves
\beq
2\Lambda_0 \ln \left( 2\Lambda_0 + \frac{1}{U}\right) -\frac{1}{U}  = t +
2\Lambda_0 \ln \left( 2\Lambda_0 + \frac{1}{\Lambda_0}\right) - \frac{1}{\Lambda_0} \,.
\eeq
Therefore, there exists $0<t_2 = t_2 (\Lambda_0) \leq 1$ such that
\be\label{Hsqbound-2}
\max_{M_t}|A| \leq \frac{3\Lambda_0}{2}\quad \text{for all } t\in [0,t_2]\,.
\ene

Then the first assertion of the Theorem, namely, \eqref{estimate-1}, is obtained from the following technical lemma
by setting $T_1 = \min \{t_1,t_2\}$.
\bl \label{t1t2}
There exists some constant $\epsilon_0 = \epsilon_0(n, |M_0|, \Lambda_0)>0$ such that if the condition
\eqref{Condition-2} is satisfied, then
\beq
t_1 \geq t_2 = t_2 (\Lambda_0)\,.
\eeq
\el
\bp{ of the Lemma ~\ref{t1t2}:} Suppose this is not the case, then we have $t_1 < t_2 \leq 1$. Then by \eqref{Hsqbound} and
\eqref{Hsqbound-2}, we deduce that at time $t =t_1$ either $h(t)$ or $\int_{M_t} H^2 d\mu$ achieves the extreme value
$2\Lambda_0$ or $\frac{1}{2\Lambda_0}$.

Now since $\{\max_{M_t}|A|, h(t)\}\leq 2\Lambda_0$ for all $t\in[0,t_1]$, integrating the equation \eqref{hdA-2} of the Corollary
\ref{hdA-1} over $M_t$, and using Hamilton's interpolation inequality for tensors (Lemma \ref{Hamiltonlemma}), we have the
{\ub} on all the higher order derivatives of $A$, which only depends on $n$ and $\Lambda_0$ (more precisely,
$\max_{t}|h(t)|, \max_{M_t}|A|$ and the initial bound on the $L^2$-norm of all the derivatives of $A$ in \eqref{Condition-1}).
In particular, for all $m \in [1,\widehat{m}]$, we have:
\be\label{gradientbound1}
\max_{M_t}|\nabla^m A|  \leq C(n, \Lambda_0) \quad \text{for }\, t\in [0,t_1]\,,
\ene
c.f. \cite[Lemma 8.3]{Hui84}.

Now we integrate the {\ee} for {\AA}, namely, the equation (iii) of the Corollary \ref{Evol-A-H} over $M_t$ for $t\in [0,t_1]$, to get
\begin{align*}
&\frac{\partial}{\partial t} \int_{M_t}|\text{\AA}|^2 \, d\mu  -  \int_{M_t}|\text{\AA}|^2 H(1-hH) \, d\mu\\
= & \int_{M_t} \left[- 2h|\nabla \text{\AA}|^2  + 2h |A|^2|\text{\AA}|^2 -2\left(\text{tr}(\text{\AA}^3)+ \frac{2}{n} H |\text{\AA}|^2\right) \right]\, d\mu\,,
\end{align*}
and therefore
\be\label{good-1}
\frac{\partial}{\partial t} \int_{M_t}|\text{\AA}|^2 \, d\mu \leq C(n, \Lambda_0)\int_{M_t}|\text{\AA}|^2 \, d\mu \quad \text{for all } t\in [0,t_1]\,,
\ene
where we have used $|H|\leq \sqrt{n}|A| \leq 2\sqrt{n}\Lambda_0$ and $|\text{tr}(\text{\AA}^3)|\leq |\text{\AA}|^3 \leq 2\Lambda_0|\text{\AA}|^2$\,.

Therefore, using \eqref{good-1} and the assumption that $\int_{M_0}|\text{\AA}|^2 \,d\mu \leq \epsilon\leq \epsilon_0$, we now get
\be\label{good2}
 \int_{M_t}|\text{\AA}|^2 \, d\mu \leq \epsilon e^{C(n, \Lambda_0)t} \leq C(n, \Lambda_0)\epsilon \quad \text{for all } t\in [0,t_1]\,,
\ene
where we abuse our notation for $C(n, \Lambda_0)$.

Now we observe from Hamilton's interpolation inequality (Lemma ~\ref{Hamiltonlemma2} with $r=1,p=q=2$):
\be\label{good-eq-6}
\int_{M_t} |\nabla \text{\AA}|^2\,d\mu \leq n \left(\int_{M_t} |\text{\AA}|^2\,d\mu\right)^{\frac{1}{2}}
\left(\int_{M_t} |\nabla^2 \text{\AA}|^2\,d\mu\right)^{\frac{1}{2}} \leq C(n,\Lambda_0)\epsilon^{\frac{1}{2}}\,,
\ene
where we used $|\nabla^2 \text{\AA}|\leq C(n)|\nabla^2 A|$ and \eqref{gradientbound1}. In fact, using \eqref{gradientbound1} and
applying Lemma \ref{Hamiltonlemma2} inductively, we have, for all $m \in [1,\widehat{m}]$,
\be
\int_{M_t} |\nabla^m \text{\AA}|^2\,d\mu \leq C(n, m, \Lambda_0)\epsilon^{\frac{1}{2m}} \quad \text{for all } t\in [0,t_1]\,.
\ene
This together with Lemma \ref{Hamiltonlemma} imply that, for all $t\in [0,t_1]$,
\beq
\int_{M_t} |\nabla^m \text{\AA}|^p\,d\mu \leq C(n, m, p, \Lambda_0)\epsilon^{\frac{1}{mp}} \quad
\text{for all }m \in [1,\widehat{m}] \text{ and } p<\infty.
\eeq
This yields, by the standard Sobolev inequality (see e.g. \cite[\S 2]{Aub98}), that for some universal constant
$\alpha \in (0,1)$, and all $m \in [1,\widehat{m}]$, and $t\in [0,t_1]$, we have
\be\label{good3}
\max_{M_t}|\nabla^m \text{\AA}| \leq C(n, m,\Lambda_0)\epsilon^{\alpha}.
\ene
In particular, using \cite[Lemma 2.2]{Hui84}, for all $t\in [0,t_1]$ we have
\be\label{Hgrad1}
\max_{M_t}|\nabla H| \leq C(n)\max_{M_t}|\nabla \text{\AA}| \leq C_1(n, \Lambda_0)\epsilon^{\alpha} \,.
\ene

Therefore, using \eqref{Hsqbound} \eqref{Hgrad1} and Topping's Theorem ~\ref{top}, we have
\begin{align}\label{hHbound02}
|1-hH|(x,t) &= \left(\int_{M_t} H^2 \, d\mu\right)^{-1} \left|\int_{M_t} H^2 \, d\mu - H(x,t)\int_{M_t} H \, d\mu \right|\notag\\
& \leq 2\Lambda_0^{-1}\text{diam}\,(M_t)\int_{M_t} |\nabla H||H|\, d\mu \notag\\
&\leq C_1(n, |M_0|, \Lambda_0)\epsilon^{\alpha}\,,
\end{align}
for all $(x,t)\in M_t$ and all $t\in [0,t_1]$. Here we abuse the notation on $C_1$ but we are allowed to choose a larger constant
$C_1$ than previously in \eqref{Hgrad1}.

Now consider the {\ee} for $H^2$, that is (see the equation (i) of Corollary \ref{Evol-A-H})
\be\label{evoH2}
\frac{\partial}{\partial t} H^2 = h\Delta H^2 - 2h|\nabla H|^2- 2(1-hH)H|A|^2\,.
\ene
Integrating this over $M_t$ we get
\be\label{int-evoH2}
\frac{\partial}{\partial t} \int_{M_t}H^2 \,d\mu = \int_{M_t} H^3(1-hH) -2h|\nabla H|^2 -2(1-hH)H|A|^2 \,d\mu\,,
\ene
where the first term on the right-hand side comes from the time derivative of $\mu_t$, i.e., the equation \eqref{mu}.

Therefore we have (using $0<\epsilon\leq 1$):
\be\label{dH1}
\left|\frac{\partial}{\partial t} \int_{M_t}H^2 \,d\mu\right| \leq C(n, \Lambda_0)|M_0|\epsilon^{\alpha} \quad \text{for all } t\in [0,t_1]\,.
\ene
Similarly, using the evolution equation for $h$, that is, the Lemma \ref{h_t},  we have
\be\label{dh1}
\left|\frac{\partial}{\partial t} h \right| \leq C(n, \Lambda_0)|M_0|\epsilon^{\alpha} \quad \text{for all } t\in [0,t_1]\,.
\ene
Integrating \eqref{dH1} and \eqref{dh1} over $[0,t_1]$ (note that $t_1<t_2\leq 1$), and choosing
$\epsilon_0 = \epsilon_0 (n, |M_0|,\Lambda_0) \geq \epsilon$ sufficiently small, we obtain
\beq
\frac{2}{3\Lambda_0} \leq \left\{ h(t_1)\,, \,\int_{M_{t_1}}H^2\,d\mu  \right\} \leq \frac{3\Lambda_0}{2}\,.
\eeq
This contradicts with the assumption that either $h(t_1)$ or $\int_{M_{t_1}} H^2 d\mu$ achieves the extreme value $2\Lambda_0$
or $\frac{1}{2\Lambda_0}$. Therefore $T_1 = t_2(\Lambda_0)$.
\ep

To see \eqref{estimate-2}, we can repeat the above argument by replacing $t_1$ by $t_2 = T_1$. Note that \eqref{good2},
\eqref{good3} and Lemma \ref{top} together yield a bound on $\max_{M_t}|\text{\AA}|$ in terms of $n, |M_0|, \Lambda_0$ and
$\epsilon$ for all $t\in [0,t_1]$. Bounds for $|\nabla H|$ and $|1-hH|$ are as in \eqref{Hgrad1} and \eqref{hHbound02}, respectively.
We will still call this bound $C_1 = C_1(n,|M_0|,\Lambda_0)$ which is chosen to be larger than the $C_1$'s in \eqref{Hgrad1} and
\eqref{hHbound02}. Now the proof is complete.
\ep
\br
Our conditions \eqref{Condition-1} appear to be necessary for the case of the {\vpmcf} in \cite{Li09} for the initial hypersurface.
Moreover, a more general form of the Kato's inequality for $|\nabla \text{\AA}|$ is probably not true, but appears to be necessary
to deploy Moser's parabolic iteration in [Pages 337, 338 and 340, \cite{Li09}]. We do not use the iteration method here to obtain
the $L^{\infty}$ bounds for $|\text{\AA}|$, $|\nabla H|$ and $|1-hH|$.
\er

\subsection{Establishing the exponential decay for geometric quantities}
Previously we have obtained a time $T_1 = T_1(\Lambda_0)$ which only depends on the initial hypersurface, and
$\epsilon_0 = \epsilon_0(n,|M_0|,\Lambda_0)$ small enough such that if the initial $L^2$ norm of {\AA} is small
(see \eqref{Condition-2}), then we have estimates \eqref{estimate-1} and \eqref{estimate-2} on time interval $[0,T_1]$. In this
subsection, we show that if on some time interval $[0,T)$, estimates similar to \eqref{estimate-1} and \eqref{estimate-2} hold,
then we can choose an $\epsilon$ small enough for the initial $L^2$ bound on {\AA}, such that $|\text{\AA}|$, $|\nabla H|$ and
$|1-hH|$ decay exponentially on this time interval $[0,T)$. More precisely,
\bt\label{exp1}
Let $M_t^n \subset \mathbb{R}^{n+1}, n\geq 2,$ be a smooth compact solution to the {\sapmcf} \eqref{eq-mcf} with initial
condition \eqref{Condition-1}. Suppose that for any $t\in [0,T), T\leq \infty$, we have
\be\label{Condition-3}
\max_{M_t}|A| \leq \Lambda_1\quad \text{and}\quad\frac{1}{\Lambda_1} \leq \left\{ h(t)\,, \,\int_{M_t}H^2\,d\mu  \right\} \leq \Lambda_1
\ene
and
\be\label{Condition-4}
\max_{M_t}\left(|\text{\AA}| + |\nabla H| + |1- hH|\right) \leq \hat{C} \epsilon^{\beta}
\ene
for some $\beta>0$. Then there exists $\epsilon_1 = \epsilon_1 (n, |M_0|, \Lambda_1, \hat{C},\beta)>0$ such that if
\be\label{Condition-5}
\int_{M_0}|\text{\AA}|^2 \,d\mu \leq \epsilon \leq \epsilon_1\,,
\ene
then for all $t\in [0,T)$ we have
\be\label{TracelessADecay1}
\max_{M_t}|\text{\AA}| \leq \left(\max_{M_0}|\text{\AA}|\right)e^{-\delta t}
\ene
and
\be\label{hHDecay1}
\max_{M_t} \left(|\text{\AA}| + |\nabla H| + | 1- hH|\right) \leq
C_2(n,|M_0|,\Lambda_1,\hat{C})\left(\max_{M_0}|\text{\AA}|\right)^{\alpha}e^{-\alpha\delta t} \,,
\ene
where $\delta = \frac{1}{4n\Lambda_1^2|M_0|} >0$, and $\alpha \in (0,1)$ is from Theorem \ref{key}.
\et
\br
To directly apply the results of Theorem ~\ref{key}, we may take $\beta = \alpha$, and $\Lambda_1 = 2\Lambda_0$,
as well as $\hat{C} = C_1$. We state the Theorem ~\ref{exp1} this way so that it can be easily adapted in later applications.
\er
\bp
Using \eqref{Condition-3}, for any $t\in [0,T)$ we have
\beq
\int_{M_t} H \, d\mu = h(t) \int_{M_t} H^2 \, d\mu \geq \frac{1}{\Lambda_1^2}\,.
\eeq
Therefore we can always find a point $x_0 \in M_t$ (which may depend on $t\in [0,T)$) such that
\be\label{Hx0}
H(x_0, t) \geq \frac{1}{\Lambda_1^2|M_t|} = \frac{1}{\Lambda_1^2 |M_0|}.
\ene
Now by \eqref{Condition-3}, \eqref{Condition-4}, the inequality $|H| \leq \sqrt{n}|A|$, and the Theorem ~\ref{top},
for fixed $t \in [0,T)$, we have
\bear
H(x,t) - H(x_0,t) &\ge& -{\max}|\nabla H| {\diam}(M_t)\\
 &\ge&  -\hat{C}\epsilon^{\beta}C(n)\int_{M_t}H^{n-1}d\mu\\
  &\ge&  -\hat{C}\epsilon^{\beta}C(n)n^{\frac {n-1}{2}}|M_0|\Lambda_1^{n-1}.
\eear
Now we can apply \eqref{Hx0}, and choose $\epsilon_1 = \epsilon_1 (n, |M_0|, \Lambda_1,\hat{C}, \beta)>0$
sufficiently small, such that if $\epsilon \leq \epsilon_1$, then for all $t\in [0,T)$, we have
\be\label{meanconvex}
{\min}_{M_t} H \geq H(x_0,t) - \hat{C}\epsilon^{\beta}C(n)n^{\frac {n-1}{2}}|M_0|\Lambda_1^{n-1}\ge \frac{1}{2\Lambda_1^2|M_0|}>0.
\ene
Moreover, we have
\bear
hH^2 &=& \frac{H^2 \int_{M_t} H d\mu}{\int_{M_t} H^2 d\mu} \\
&\leq& \frac{H^2 \int_{M_t} \{H(x_0,t) + {\max}|\nabla H| {\diam}(M_t)\}d\mu}{\int_{M_t}\{H(x_0,t) - {\max}|\nabla H| {\diam}(M_t)\}^2 d\mu}\\
&\leq& \frac{\{H(x_0,t) + {\max}|\nabla H| {\diam}(M_t)\}^3}{\{H(x_0,t) - {\max}|\nabla H| {\diam}(M_t)\}^2}.
\eear
We can then apply the estimate on $|\nabla H|$, namely \eqref{Condition-4}, we can choose $\epsilon_1$ small enough
such that from above, we have
\beq
hH^2 \leq \frac32 H(x_0,t) - 2\max(|\nabla H|)\diam(M_t).
\eeq
This implies, by choosing a possibly smaller $\epsilon_1$, we have
\begin{align}\label{good1}
\max_{M_t}\left(\frac{2}{n} h H^2 - \frac{4}{n}H\right) &= \max_{M_t}\frac{2}{n} \left( hH^2-2 H\right) \notag \\
&\leq \frac2n\left(\frac32 H(x_0,t) - 2H(x_0,t) \right)\notag\\
&\leq - \frac{1}{n\Lambda_1^2|M_0|}\,.
\end{align}
Here in the last step, we applied \eqref{Hx0}.

To derive exponential decay for $|\text{\AA}|^2$, we recall its {\ee}, namely, (iii) of Corollary \ref{Evol-A-H}, we have
\bear
\frac{\partial}{\partial t} |\text{\AA}|^2 &=& h\Delta |\text{\AA}|^2 - 2h|\nabla \text{\AA}|^2 +
2h|A|^2|\text{\AA}|^2-2\left(\text{tr}(\text{\AA}^3)+ \frac{2}{n} H |\text{\AA}|^2\right)\\
&\leq& h\Delta |\text{\AA}|^2 + \left(2h(|\text{\AA}|^2 + \frac{1}{n}H^2) + 2|\text{\AA}| - \frac{4}{n}H\right)|\text{\AA}|^2 \\
&=& h\Delta |\text{\AA}|^2 + (2h|\text{\AA}|^2 + 2|\text{\AA}|)|\text{\AA}|^2  +\left(\frac{2}{n}hH^2 - \frac{4}{n}H\right)|\text{\AA}|^2
\eear
Now we use the inequality $|\text{\AA}| \le \hat{C}\epsilon^\beta$ (\eqref{Condition-4}), and the inequality \eqref{good1}, by
choosing a possibly smaller $\epsilon_1$, we have for $\epsilon < \epsilon_1$,
\be
\frac{\partial}{\partial t} |\text{\AA}|^2 \leq h\Delta |\text{\AA}|^2 - \frac{1}{2n\Lambda_1^2|M_0|}|\text{\AA}|^2
= h\Delta |\text{\AA}|^2 - 2\delta|\text{\AA}|^2,
\ene
where $\delta = \frac{1}{4n\Lambda_1^2|M_0|}$.

Therefore the exponential decay of $|\text{\AA}|$, namely the estimate \eqref{TracelessADecay1} now follows from the {\maxp}
(Theorem \ref{maxPrin}).

Finally, once we obtain \eqref{TracelessADecay1}, we can prove \eqref{hHDecay1} exactly following the argument in the
proof of the Theorem ~\ref{key} (see \eqref{good-eq-6}--\eqref{hHbound02}).
\end{proof}
\section{Proof of Theorem \ref{MainThm}: continued}
We now assemble estimates obtained from last section to complete the proof of the main theorem: in \S 4.1, we prove the
long-time existence of the flow \eqref{eq-mcf} by establishing the uniform upper bound for $|A|$; in \S 4.2, we show the
exponential convergence of the flow.

\subsection{Extending the time interval}
In the previous section, we obtain the exponential decay for $|\text{\AA}|$, $|\nabla H|$ and $|1-hH|$ on some time interval.
We will next show that this implies a uniform bound on the function $h$, which consequently yields the {\ub} on $|A|$. We
will state the theorem in a more general form in order for later application.
\bt\label{exp2}
Let $M_t^n \subset \mathbb{R}^{n+1}, n\geq 2,$ be a smooth compact solution to the {\sapmcf} \eqref{eq-mcf} with
initial conditions \eqref{Condition-1} and \eqref{Condition-2}. Suppose that for any $t\in [0,T), T\leq \infty$ we have
\be\label{Condition-6}
\max_{M_t}|A| \leq \Lambda_2\quad \text{and}\quad \frac{1}{\Lambda_2} \leq \left\{h(t)\,, \,\int_{M_t}H^2\,d\mu  \right\} \leq \Lambda_2
\ene
and
\be\label{Condition-7}
\max_{M_t}\left(|\text{\AA}| + |\nabla H| + | 1- hH|\right)  \leq \tilde{C} \epsilon^{\beta} e^{-\alpha\delta t}\,,
\ene
for some $\beta, \delta >0 $ and $\alpha>0$, where $\alpha$ is the same as in the Theorem \ref{key}, and $\beta, \delta$
are the same as in the Theorem \ref{exp1}.

Then we have the following uniform estimate for all $t\in [0,T)$:
\be\label{ub-h0}
\int_{M_0}H^2 \,d\mu - b_0\delta^{-1}\epsilon^{\beta} \leq \int_{M_{t}}H^2 \,d\mu \leq \int_{M_0}H^2 \,d\mu
+ b_0\delta^{-1}\epsilon^{\beta} \,,
\ene
where $b_0 = b_0(n,|M_0|,\Lambda_2, \tilde{C}) = 2(n^{\frac32}\Lambda_2^2 + \tilde{C})\Lambda_2|M_0|\tilde{C}$ and
\be\label{ub-h1}
h(0) - b_1 \delta^{-1}\epsilon^{\beta} \leq h(t) \leq h(0) + b_1 \delta^{-1}\epsilon^{\beta}\,,
\ene
where
$b_1 = b_1(n,|M_0|,\Lambda_2, \tilde{C}) =(2n^{\frac12}\Lambda_2^2 + 2n \tilde{C} +1)\Lambda_2^3|M_0|\tilde{C} $.

Moreover, there exists $\epsilon_2 = \epsilon_2 (n, |M_0|, \Lambda_2, \tilde{C},\beta,\delta)>0$ such that if
$\epsilon \leq \epsilon_2$ then for any $t\in[0,T)$
\be\label{ub-h101}
\max_{M_t} |A| \leq  \Lambda_0^4,
\ene
and
\be\label{ub-h102}
\frac{1}{\Lambda_0^4}\leq \left\{ h(t)\,, \,\int_{M_t}H^2\,d\mu\,,\,\int_{M_t}|\nabla^m A|^2\,d\mu  \right\}\leq  \Lambda_0^4,
\end{equation}
for all $m \in [1,\widehat{m}]$, where $\widehat{m}$ and $\Lambda_0$ are from the Theorem \ref{key}.
\et

\bp
Without loss of generality, we assume $\epsilon \leq 1$. We start by recalling the following integral \eqref{int-evoH2}:
\beq
\frac{\partial}{\partial t} \int_{M_t}H^2 \,d\mu = \int_{M_t} H^3(1-hH) -2h|\nabla H|^2 -2(1-hH)H|A|^2 \,d\mu.
\eeq
Using the assumptions \eqref{Condition-6} and \eqref{Condition-7}, we can estimate this integral term by term as follows:
\beq
|H^3(1-hH)| \leq n^{\frac32}\Lambda_2^3\tilde{C}\epsilon^\beta e^{-\delta t},
\eeq
and
\beq
2h|\nabla H|^2 \leq 2\Lambda_2\tilde{C}^2\epsilon^{2\beta}e^{-2\delta t},
\eeq
and
\beq
2H|1-hH||A|^2 \leq 2n^{\frac12}\Lambda_2^3\tilde{C}\epsilon^\beta e^{-\delta t}.
\eeq
Putting these estimates together, and abusing our notation for $\tilde{C}$ and $\delta$, we obtain
\be
\left|\frac{\partial}{\partial t} \int_{M_t}H^2 \,d\mu\right| \leq \tilde{C}\epsilon^{\beta} e^{-\delta t}\quad \text{for all } t\in [0,T).
\ene
Integrating this over $[0,t]$ for any $t\leq T$ we obtain the estimate \eqref{ub-h0}.

In order to show the estimate \eqref{ub-h1}, we use the {\ee} of $h$, namely, the Lemma ~\ref{h_t}:
\beq
\frac{\partial}{\partial t} h = \frac{\int_{M_t}[-(1-2hH)(1-hH)|A|^2 + H^2(1-hH)^2 + 2h^2|\nabla H|^2\, ]d\mu}{\int_{M_t} H^2 \, d\mu}.
\eeq
We again estimate it term by term, under the assumptions \eqref{Condition-6}, \eqref{Condition-7}. Up to abuse of the notation for
$\tilde{C}= \tilde{C}(n,|M_0|,\Lambda_2)$ and $\delta$,  assuming again that $\epsilon \leq 1$, we have for any $t\in [0,T)$:
\be\label{ub-ht}
\left|\frac{\partial}{\partial t} h\right| \leq \tilde{C}\epsilon^{\beta} e^{-\delta t}.
\ene
Integrating this over $[0,t]$ for any $t\leq T$ we get \eqref{ub-h1}.

We then use
\beq
|H| \leq \frac{|1-hH|}{h} + \frac1h,
\eeq
and we choose $\epsilon_2 = \epsilon_2 (n, |M_0|, \Lambda_2, \tilde{C},\beta,\delta)>0$  sufficiently small, in view of the initial condition
\eqref{Condition-1} ($\frac{1}{\Lambda_0} \leq h(0) \leq \Lambda_0$), and \eqref{Condition-7}--\eqref{ub-h1}, we get, for any $t\in [0,T)$,
\bear
\max_{M_t}|H| &\leq& \frac{|1-hH|}{h} + \frac1h \\
&\leq& \frac{5}{4h} \\
&\leq& \frac{5}{4(h(0) - \frac{1}{4\Lambda_0})} \\
&\leq& \frac{5\Lambda_0}{3}.
\eear
and thus
\be
\max_{M_t}|A| \leq \max_{M_t}\left(\sqrt{|\text{\AA}|^2 + \frac1n |H|^2}\right) \leq 2\Lambda_0 \leq \Lambda_0^4,
\ene
where we can choose $\Lambda_0$ large enough for the last inequality. This proves \eqref{ub-h101}.

The bound on $h(t)$ and $\int_{M_t}H^2\,d\mu$ in \eqref{ub-h102} follows immediately from \eqref{Condition-1}, \eqref{ub-h0}
and \eqref{ub-h1}. We are left to estimate the integral $\int_{M_t}|\nabla^m A|^2\,d\mu$.

Now the Lemma \ref{interpolation01} yields for all $m \in [1,\widehat{m}]$,
\begin{eqnarray}\label{int-nabla}
\frac{d}{dt} \int_{M_t}|\nabla^m A|^2\,d\mu &\leq& C(n,m,|h|)\max_{M_t}\{|A|^2 + |A|\} \int_{M_t}|\nabla^m A|^2\,d\mu \nonumber\\
&\leq& C'(n,m,|h|)\Lambda_0^2 \int_{M_t}|\nabla^m A|^2\,d\mu.
\end{eqnarray}
Since we are allowed to assume $\Lambda_0$ in Theorem ~\ref{key} is sufficiently large so that
\beq
\Lambda_0 \geq \max \{ C'(n, m, |h|), 100\}\,,
\eeq
then we obtain from the inequality \eqref{int-nabla} that
\beq
\int_{M_t}|\nabla^m A|^2\,d\mu \leq C'(n, m, |h|) \Lambda_0^2 \int_{M_0}|\nabla^m A|^2\,d\mu \leq \Lambda_0^4\,.
\eeq
Now our proof is complete.
\ep
We now complete the proof for long-time existence of the flow by the following extension theorem:
\bt\label{ext}
Let $M_t^n \subset \mathbb{R}^{n+1}, n\geq 2,$ be a smooth compact solution to the {\sapmcf} \eqref{eq-mcf}
with initial conditions \eqref{Condition-1} and \eqref{Condition-2}. Suppose that for any $t\in [0,T], T<\infty$ and
for all $m \in [1, \widehat{m}]$ we have
\be\label{wholebound-1}
 \max_{M_t}|A| \leq  \Lambda_0^4 \quad \text{and}\quad \frac{1}{\Lambda_0^4} \leq
 \left\{ h(t)\,, \,\int_{M_t}H^2\,d\mu\,,\, \int_{M_t}|\nabla^m A|^2\,d\mu \right\}\leq  \Lambda_0^4
\ene
and
\be\label{wholebound-2}
\max_{M_t}\left(|\text{\AA}| + |\nabla H| + | 1- hH|\right) \leq
C^{\ast} \epsilon^\frac{\alpha^2}{2} e^{-\alpha\delta t} \leq C^{\ast} \epsilon^\frac{\alpha^2}{2}\,.
\ene
Here $0 < \alpha < 1$ is the universal constant from the Theorem ~\ref{key} and
$\delta = \frac{1}{16n\Lambda_0^8|M_0|} >0$. Then there exists some
$\epsilon_3 = \epsilon_3(n, |M_0|, \Lambda_0, \alpha, C^\ast)>0$ and $T_2 = T_2(\Lambda_0) >0$ such that if
\be
\int_{M_0}|\text{\AA}|^2 \,d\mu \leq \epsilon \leq \epsilon_3\,,
\ene
then \eqref{wholebound-1} and \eqref{wholebound-2} hold for all $t\in [0,T+T_2]$.
\et
\bp
We start by applying the Theorem ~\ref{key} while setting the initial time to be $t = T$. Then there exists
$\epsilon_4 := \epsilon_0 (n, |M_0|, \Lambda_0^4) >0$ and $T_2:=  T_1(\Lambda_0^4)>0$ such that if
\beq
\int_{M_0}|\text{\AA}|^2 \,d\mu \leq \epsilon \leq \epsilon_3\,,
\eeq
then for all $t\in [T,T+T_2]$ we have
\be\label{good4}
\max_{M_t}|A| \leq  2\Lambda_0^4\quad \text{and}\quad
\frac{1}{2\Lambda_0^4} \leq \left\{ h(t)\,, \,\int_{M_t}H^2\,d\mu  \right\}\leq  2\Lambda_0^4,
\ene
and for some $\alpha \in (0,1)$,
\be\label{good5}
\max_{M_t}\left(|\text{\AA}| + |\nabla H| + |1- hH|\right) \leq   C_1(n,|M_0|, \Lambda_0^4) \epsilon^{\alpha}\,.
\ene

Now we can choose $\epsilon_5 = \epsilon_5 (n, |M_0|, \Lambda_0, \alpha) >0$ sufficiently small such that
for all $\epsilon < \epsilon_5$, we have
\beq
C_1(n,|M_0|, \Lambda_0^4) \epsilon^{\alpha -\frac{\alpha^2}{2}} \leq C^{\ast}\,
\eeq
and therefore for all $t\in [0,T+T_2]$, we have
\be\label{good7}
\max_{M_t}\left(|\text{\AA}| + |\nabla H| + | 1- hH|\right)  \leq C^{\ast} \epsilon^\frac{\alpha^2}{2}\,.
\ene

We are now in position to apply the Theorem ~\ref{exp1} as follows. Given \eqref{wholebound-1} for $t\in [0,T]$,
\eqref{good4} for $t\in [T,T+T_2]$ and \eqref{good7} for $t\in [0,T+T_2]$, we apply the Theorem ~\ref{exp1} on the
time interval $[0, T+T_2]$ with $\Lambda_1 = 2\Lambda_0^4, \hat{C} = C^\ast$ and $\beta = \frac{\alpha^2}{2}$ to
find that there exists some $\epsilon_6:= \epsilon_1 (n, |M_0|,\Lambda_0, C^\ast,\alpha)>0$ sufficiently small,
so that if $\epsilon \leq \epsilon_6$, then for all $t\in [0,T+T_2]$, we have
\begin{eqnarray}\label{exp3}
\max_{M_t} \left(|\text{\AA}| + |\nabla H| + | 1- hH|\right) &\leq&
C_2(n,|M_0|,2\Lambda_0^4,C^\ast)\left(\max_{M_0}|\text{\AA}|\right)^{\alpha}e^{-\alpha\delta t} \nonumber \\
&\leq& C_2(n,|M_0|,2\Lambda_0^4,C^\ast) \left(C_1 (n, |M_0|, \Lambda_0)\right)^{\alpha} \epsilon^{{\alpha}^2}e^{-\alpha\delta t}\,,
\end{eqnarray}
where $\delta = \frac{1}{16n \Lambda_0^8|M_0|} >0$. Here we have also used the estimate \eqref{estimate-2} at $t=0$.

We can then proceed to apply the Theorem ~\ref{exp2}. To do so, we first choose some
$\epsilon_7 = \epsilon_7 (n, |M_0|, \Lambda_0, \alpha, C^{\ast}) >0$ small enough, so that
\be\label{good8}
 C_2(n,|M_0|,2\Lambda_0^4,C^\ast) \left(C_1 (n, |M_0|, \Lambda_0)\right)^{\alpha} \epsilon^{\frac{{\alpha}^2}{2}} \leq C^{\ast}.
\ene
This allows us to rewrite the estimate \eqref{exp3} to:
\be\label{exp4}
\max_{M_t} \left(|\text{\AA}| + |\nabla H| + | 1- hH|\right) \leq C^{\ast}\epsilon^{\frac{{\alpha}^2}{2}}e^{-\alpha\delta t}.
\ene

Comparing with \eqref{Condition-6} and \eqref{Condition-7}, we can then apply the Theorem ~\ref{exp2} to the time interval
$[0,T+T_2]$ with $\Lambda_2 = 2\Lambda_0^4$, $\beta = \frac{\alpha^2}{2}$, $\delta = \frac{1}{16n \Lambda_0^8|M_0|}$,
and $\tilde{C} = C^{\ast}$, so that we can choose
\beq
\epsilon_8 = \epsilon_8 (n, |M_0|,\Lambda_0,\alpha, C^\ast) :=
\epsilon_2 \left(n,|M_0|, 2\Lambda_0^4, C^{\ast}, \frac{\alpha^2}{2}, \frac{1}{16n \Lambda_0^8|M_0|}\right)
\eeq
and if $\epsilon \leq \epsilon_8$ then we have \eqref{wholebound-1}.

We complete the proof by setting
\beq
\epsilon_3 = \epsilon_3(n, |M_0|, \Lambda_0, \alpha, C^\ast) = \min \left\{ \epsilon_4, \epsilon_5, \epsilon_6 , \epsilon_7, \epsilon_8 \right\} > 0 \,.
\eeq
\ep
\subsection{Completion of the proof}
We now complete the proof of our main theorem:
\begin{proof}(of Theorem \ref{MainThm}) Suppose that the initial condition \eqref{Condition-1} is satisfied for some
$\Lambda_0 \geq 100$ sufficiently large. Then by the Theorem ~\ref{key}, we first choose
$\epsilon_0 = \epsilon_0(n, |M_0|, \Lambda_0)>0$ and $T_1 = T_1 (\Lambda_0) \in (0,1]$, such that if
$\epsilon \leq \epsilon_0$, then estimates \eqref{estimate-1} and \eqref{estimate-2} hold for all $t\in [0,T_1]$.

Now we can apply the Theorem ~\ref{exp1} to the interval $[0,T_1]$ with
$\Lambda_1 = 2\Lambda_0$, $\hat{C} = C_1 (n, |M_0|, \Lambda_0)$, and $\beta = \alpha$, for some
\beq
\epsilon_9 = \epsilon_9 (n, |M_0|, \Lambda_0, \alpha) := \epsilon_1\left(n, |M_0|, 2\Lambda_0, C_1 (n, |M_0|, \Lambda_0), \alpha\right)>0
\eeq
sufficiently small such that if $\epsilon\leq \epsilon_9$, then for all $t\in [0,T_1]$, we have
\begin{align}\label{last-ineq}
&\max_{M_t} \left(|\text{\AA}| + |\nabla H| + | 1- hH|\right) \notag \\
\leq &\, C_2\left(n, |M_0|, 2\Lambda_0, C_1(n, |M_0|, \Lambda_0)\right)\left(\max_{M_0}|\text{\AA}|\right)^{\alpha}e^{-\alpha\gamma t}\notag \\
\leq & \,C_2\left(n, |M_0|, 2\Lambda_0, C_1(n, |M_0|, \Lambda_0)\right)\left(C_1 (n, |M_0|, \Lambda_0)\right)^{\alpha} \epsilon^{\frac{{\alpha}^2}{2}} e^{-\alpha\delta t}\,,
\end{align}
where $\gamma = \frac{1}{16n\Lambda_0^2|M_0|} \geq \delta  = \frac{1}{16n\Lambda_0^8|M_0|}>0$,
and we have again used the estimate \eqref{estimate-2} at $t=0$.

Let $C^\ast = C_2\left(n, |M_0|, 2\Lambda_0, C_1 (n, |M_0|, \Lambda_0)\right)\left(C_1 (n, |M_0|, \Lambda_0)\right)^{\alpha}$, the above
inequality \eqref{last-ineq} becomes
\be
\max_{M_t} \left(|\text{\AA}| + |\nabla H| + | 1- hH|\right) \leq C^\ast\epsilon^{\frac{{\alpha}^2}{2}} e^{-\alpha\delta t}.
\ene

This allows us to apply the Theorem ~\ref{ext}. We see that if we choose (note that $0 < \alpha < 1$ is some universal constant)
\beq
\epsilon \leq \epsilon_{10} = \epsilon_{10} (n, |M_0|, \Lambda_0) := \min \{ \epsilon_9 , \epsilon_3(n, |M_0|, \Lambda_0, \alpha, C^\ast)\}
\eeq
then the flow \eqref{eq-mcf} exists for all time and converges exponentially to a round sphere.
\ep

\bibliographystyle{amsalpha}
\bibliography{ref-SAP}
\end{document}